\documentclass[12pt]{amsart}
\usepackage{graphicx}
\usepackage{amssymb}
\usepackage{amsfonts}
\usepackage{amsmath}
\usepackage{array}
\usepackage{rotating}
\usepackage{bm}
\usepackage[matrix,frame]{xypic}
\usepackage{MnSymbol}
\newtheorem{example}{Example}[section]

\newtheorem{theorem}[example]{Theorem}

\newtheorem{proposition}[example]{Proposition}

\newtheorem{lemma}[example]{Lemma}

\def\Proof{\noindent \it Proof -- \rm}
\def\qed{\hspace{3.5mm} \hfill \vbox{\hrule height 3pt depth 2 pt width 2mm}
\bigskip}

\def\PW{{\rm PW}}   

\def\stuffle{\uplus}
\def\MS{{\bf MS}}
\def\XX{{\bf X}}

\def\MQSym{{\bf MQSym}}

\def\WQSym{{\bf WQSym}}
\def\Sym{{\bf Sym}}
\def\NCSF{{\bf Sym}}
\def\QSym{{\it QSym}}

\def\ev{{\rm Ev}}

\def\<{\langle}
\def\>{\rangle}

\def\K{\operatorname{\mathbb K}}

\def\N{\operatorname{\mathbb N}}

\def\M{{\bf M}}

\def\SG{{\mathfrak S}}

\def\A{{\sf A}}

\def\shuff#1#2{\mathbin{
\hbox{\vbox{ \hbox{\vrule \hskip#2 \vrule height#1 width 0pt
}%
\hrule}%
\vbox{ \hbox{\vrule \hskip#2 \vrule height#1 width 0pt
\vrule }%
\hrule}%
}}}

\def\shuf{{\mathchoice{\shuff{7pt}{3.5pt}}%
{\shuff{6pt}{3pt}}%
{\shuff{4pt}{2pt}}%
{\shuff{3pt}{1.5pt}}}}%

\def\bshuf{\,\hat\shuf\,}

\def\npn{{\bf n}}

\def\A{{\bm A}}

\def\X{{\mathbb X}}

\def\AA{{\mathcal A}} 

\def\pack{{\rm pack}}

\title[]%
{A two-parameter deformation of the quasi-shuffle\\ and new bases of quasi-symmetric functions}
\author[O. Bouillot, J.-C.~Novelli and J.-Y.~Thibon]%
{Olivier Bouillot, Jean-Christophe Novelli and Jean-Yves Thibon}

\address[]{[Bouillot, Novelli, Thibon] Laboratoire d'informatique Gaspard-Monge\\
Universit\'e Gustave Eiffel \\
5, Boulevard Descartes \\ Champs-sur-Marne \\
77454 Marne-la-Vall\'ee cedex 2 \\
France}
\email[Olivier Bouillot]{olivier.bouillot@univ-eiffel.fr}
\email[Jean-Christophe Novelli]{jean-christophe.novelli@univ-eiffel.fr}
\email[Jean-Yves Thibon]{jean-yves.thibon@univ-eiffel.fr} 
\date{\today}


\date{\today}

\begin{document}

\begin{abstract}
	We define a two-parameter deformation of the quasi-shuffle by means of the formal group law associated with the
	exponential generating function of the homogeneous Eulerian polynomials, and construct bases of $QSym$ and $\WQSym$ 
	whose product rule is given by this operation. 
\end{abstract}

\maketitle

\section{Introduction}

The Hopf algebra of quasi-symmetric functions \cite{Ge} is historically the first example
of a combinatorial Hopf algebra extending that of symmetric functions.
While there is no general agreement on what should be the formal definition of a combinatorial Hopf algebra,
its is clear that to be considered as combinatorial, such an algebra should, in addition to  combinatorial product and coproduct
rules, also have at least two bases related by a combinatorial rule.

Quasi-symmetric functions entered the scene with these requirements, having from the beginning two bases
$M_I$ (monomial) and $F_I$ (fundamental). The product rule for the $M_I$ is the so-called
quasi-shuffle (apparently first described in \cite{TU}), and its coproduct is given by deconcatenation. 

At the time, the algebra of symmetric functions had already a lot of known bases
(elementary, complete, Schur, power-sum, Hall-Littlewood, Jack, Macdonald etc., see \cite{Mac}), and the
subsequent realization of the dual of $QSym$ as noncommutative symmetric functions
quickly led to the introduction of a lot of new bases on both sides (see {\it e.g.},\cite{NCSF2,HNTT,NTT,LNT,NTW, SchaThi,NTbinbsh,NTToum,NTLag2}
and the survey  \cite{Mason} for a sample).

Originally motivated by the study of posets, quasi-symmetric functions were next related
to the combinatorics of descents in permutations and to the Solomon descent algebra,
then to the 0-Hecke algebras and certain degeneracies of quantum groups \cite{NCSF4}.
More recently, they played a crucial role in the theory of Macdonald polynomials,
via the quasi-symmetric expansions of certain LLT polynomials \cite{HHL}.

Another source of interest in quasi-symmetric functions stems in the theory of
multiple zeta values, which are the specializations $x_n= \frac1n$ of
the monomial functions $M_I$. Thus, the multiple zeta values satisfy a product formula given by the
quasi-shuffle, but their representation as iterated integrals naturally labelled
by binary words leads to another expression of the product, given by the ordinary
shuffle of binary words. The equality between both expressions 
encodes  algebraic relations satisfied by the multiple zeta values \cite{Zag}.

A natural question is therefore: does there exist bases of $QSym$ naturally labelled by binary words
which multiply by the ordinary shuffle of binary words? This question has been answered affirmatively
in \cite{NTbinbsh}, where a general method for constructing such bases has been presented, and a special
case investigated more in depth.

Recently, a variant of the quasi-shuffle called the block shuffle has been introduced
by Hirose and Sato, also in relation with the multiple zeta values. It has been proved by Keilthy \cite{Kei}
that the block shuffle algebra is isomorphic to the ordinary shuffle algebra, providing an isomorphism
analogous to Hoffman's exponential \cite{Hof}, where the exponential is replaced by the hyperbolic tangent.

Actually, Hoffman's classical proof of the isomorphism between shuffle and quasi-shuffle algebras
is equivalent to the fact that $QSym$ admits bases whose product rule is given by the ordinary
shuffle of compositions\footnote{This proves the result for the quasi-shuffle algebra over the semigroup
of nonnegative integers, but the argument is easily extended to quasi-symmetric functions with exponents
in an appropriate additive semigroup \cite{NPT}}. The existence of such bases comes itself from the fact that
$\NCSF$ is the free graded associative algebra over a sequence of primitive elements.

Once more, a natural question is now to construct bases of quasi-symmetric functions whose
product rule is given by the block shuffle. The existence of such bases is implied by
the result of \cite{Kei}. We shall reprove this result in a simpler way, starting on the
dual side with a sequence of primitive generators, and investigate the bases obtained
by a particular choice of these generators, corresponding to the Solomon idempotents
of the descent algebra.

For a word $w=a_1\cdots a_n$ over the alphabet of positive integers, and an integer $x$, define
$\zeta_x(w)=(a_1+x)a_2\cdots a_n$, where $a_1+x$ is the sum of integers.
We shall  see that the ordinary shuffle, the quasi-shuffle and the block shuffle are
special cases of a two-parameter family of associative products defined by the recurrence
\begin{equation}
	au\star bv = a(u\star bv)+b(au\star v)+\alpha[a+b]\cdot (u\star v) +\beta\zeta_{a+b}(u\star v)
\end{equation}
associated with the formal group law
\begin{equation}
	F(x,y)=\frac{x+y+\alpha xy}{1-\beta xy}.
\end{equation}
This is proved by introducing a basis of $QSym$ satifying this product rule.
This basis is then lifted to $\WQSym$.

\subsection*{Acknowledgements} This research has been partially supported by the
program CARPLO (ANR-20-CE40-0007) of the Agence Nationale de la Recheche.

\bigskip
We shall employ the notations of \cite{NCSF1} and \cite{NTamaz}. For the convenience of the reader, the necessary background is collected in an appendix (Section \ref{sec:app}).

\goodbreak
\section{The block shuffle}

The block shuffle of words over the alphabet of positive integers is defined by

\begin{equation}
	w\bshuf\epsilon=\epsilon\bshuf w = w
\end{equation}
where $\epsilon$ is the empty word, and for letters $a,b$

\begin{equation}
	au\bshuf bv = a(u\bshuf bv)+b(au\bshuf v)-\zeta_{a+b}(u\bshuf v)
\end{equation}
where $\zeta_a(bv)=(a+b)v$ ($a+b$ is the sum of the integers) and $\zeta_a(\epsilon)=0$.

Suppose that $X_I$ is a basis of $QSym$ satifying $X_IX_J=X(I\bshuf J)$
(where $X$ is regarded as a linear map on the vector space spanned by compositions, and $X(I):=X_I$).
Let $Y_I$ be the dual basis of $X_I$. From the product formula of the $X_I$, we can calculate
the coproduct of $Y_I$. In particular, to get  $\Delta Y_n$, we see that $X_n$ can appear in $X_IX_J$ if and only if
$|\ell(I)-\ell(J)|=1$, and that its coefficient is a sign $(-1)^r$, where $r=\frac12(\ell(I)+\ell(J)-1)$.
Hence,
\begin{equation}\label{eq:DY}
	\Delta Y_n = \sum_{|I|+|J|=n\atop |\ell(I)-\ell(J)|=1}(-1)^{(\ell(I)+\ell(J)-1)/2}Y_I\otimes Y_J.
\end{equation}
Assuming that $Y_I=Y_{i_1}\cdots Y_{i_r}=:Y^I$ is a multiplicative basis, and setting $Y=\sum_{n\ge 1}Y_n$, this can be rewritten as
\begin{equation}
\Delta Y=\frac{Y\otimes 1+1\otimes Y}{1\otimes 1+Y\otimes Y}.
\end{equation} 
Starting with a sequence $\Pi_n$ of primitive generators of $\Sym$, and setting
\begin{equation}
	Y=\sum_{n\ge 1} Y_n=\tanh \Pi\quad \text{where}\  \Pi=\sum_{n\ge 1}\Pi_n
\end{equation}
we have from the addition formula
\begin{equation}
	\tanh(a+b)=\frac{\tanh a+\tanh b}{1+\tanh a\tanh b}
\end{equation}

\begin{equation}
	\Delta Y = Y(A+B)= \tanh(\Pi(A)+\Pi(B)) = \frac{Y\otimes 1+1\otimes Y}{1\otimes 1+Y\otimes Y}
\end{equation}
which is the generating series of \eqref{eq:DY}. Thus, defining $Y_n$ in this way does indeed yield a basis of $\NCSF$
whose dual basis $X_I$ mutiplies by the block shuffle.

\section{A special choice}\label{sec:spec}

Let us now choose $\Pi=\varphi$, where

\begin{equation}
	\varphi=\log \sigma_1 = \sum_{n\ge 1}\varphi_n = \sum_{n\ge 1}\frac{\Phi_n}{n}
\end{equation}
so that (recall that $\lambda_{-1}=\sigma_1^{-1}$)
\begin{equation}
	Y =\tanh \varphi =\frac{e^\varphi-e^{-\varphi}}{e^\varphi+e^{-\varphi}}=
 \frac{\sigma_1-\lambda_{-1}}{\sigma_1+\lambda_{-1}}=\frac{\sigma_1(2A)-1}{2}\left(1+\frac{\sigma_1(2A)-1}{2}\right)^{-1}
\end{equation}
whence
\begin{align}
	Y_n &= \sum_{I\vDash n}\frac{(-1)^{\ell(I)-1}}{2^{\ell(I)}}S^I(2A)\\
	&= \sum_{I\vDash n}\frac{(-1)^{\ell(I)-1}}{2^{\ell(I)}}\sum_{J\le I}R_J(2A)\\
	&= \sum_{J\vDash n}\sum_{I\ge J}\frac{(-1)^{\ell(I)-1}}{2^{\ell(I)}}R_J(2A)\\
	&=\frac1{2^n} \sum_{J\vDash n}(-1)^{\ell(J)-1}R_J(2A).
\end{align}
Thus,

\begin{equation}
	Y_n = \frac{-1}{2^n}\AA_n(-1;2A),
\end{equation}
where 
\begin{equation}
	\AA_n(t;A) = \sum_{I\vDash n}t^{\ell(I)}R_I(A)=\sum_{k=1}^n t^k\A(n,k)
\end{equation}
are the noncommutative Eulerian polynomials of \cite{NCSF1,NTamaz}. Therefore,
with the notation of \cite{NCSF1,NTamaz}
\begin{equation}
	Y_n=\frac1{2^n}\sum_{k=1}^n(-1)^{k-1}\A(n,k)*S_n^{[2]},
\end{equation}
and using \cite[(16)]{NTamaz}, we can compute
\begin{align}
	\A(n,k)*S_n^{[2]}& =\sum_{i=0}^k(-1)^i{n+1\choose i}S_n^{[k]}*S_n^{[2]}\\
	&= \sum_{i=0}^k(-1)^i{n+1\choose i}\sum_{j=0}^{2k-2i}{n+j\choose j}\A(n,2k-2i-j)\\
	&=\sum_j\sum_{i=0}^k(-1)^i {n+j\choose j}{n+1\choose i}\A(n,2k-2i-j)\\
	&=\sum_r \left( \sum_{i=0}^k(-1)^i{n+1\choose i}{n+2k-2i-r\choose n}\right)\A(n,r).
\end{align}

\begin{lemma}
	The coefficient of $\A(n,r)$ in the above sum is

\begin{equation} \label{eq-fnkr}
	F(n,k,r):= \sum_{i=0}^k(-1)^i{n+1\choose i}{n+2k-2i-r\choose n} = {n+1\choose 2k-r}.
\end{equation}
\end{lemma}

\Proof 
One easily checks that
\begin{equation}
F(n,k,r) = F(n,k,r+1) + F(n-1,k,r) - F(n-1,k,r+2).
\end{equation}
This equation is also obviously satisfied by the right-hand side
of~\eqref{eq-fnkr}. Moreover, the initial conditions are
\begin{itemize}
\item $r=2k$, where Equation~\eqref{eq-fnkr} reduces to $1=1$,
\item $r=2k-1$, where Equation~\eqref{eq-fnkr} reduces to $n+1=n+1$,
\item $n=0$, where Equation~\eqref{eq-fnkr} reduces to
$\binom{2k-r}{0} - \binom{2k-r-2}{0} = \binom{1}{2k-r}$ which is true for all
values of $2k-r$.
\end{itemize}
\qed

Finally, we arrive at the following expression:
\begin{theorem}
On the ribbon basis,
\begin{equation}
	Y_n=\frac1{2^n}\sum_{r=1}^n \Re \left(i^r(1+i)^{n+1}\right)\A(n,r),
\end{equation}
where $\Re$ denotes the real part,
or, equivalently
\begin{equation}
	Y_n=\left(\frac1{\sqrt{2}}\right)^{n-1}\sum_{I\vDash n}\cos (n+1+2\ell(I))\frac{\pi}{4}\cdot R_I.
\end{equation}
\end{theorem}

\Proof The coefficient of $\A(n,r)$ in $Y_n$ is $2^{-n}a(n,r)$, where
\begin{equation}
	a(n,r) = \sum_{k=1}^n(-1)^k{n+1\choose 2k-r}=\Re\sum_{p\ge 0}i^{p+r}{n+1\choose p}=\Re\left( i^r(1+i)^{n+1}\right).
\end{equation}
\qed

For example,
\begin{align}
	Y_2&=-\frac{1}{2}R_{1 1} + \frac{1}{2}R_{2}\\
	Y_3&=-\frac{1}{2}R_{1 2} - \frac{1}{2}R_{2 1}\\
	Y_4&=\frac{1}{4}R_{1 1 1 1} + \frac{1}{4}R_{1 1 2} + \frac{1}{4}R_{1 2 1} - \frac{1}{4}R_{1 3} + \frac{1}{4}R_{2 1 1} - \frac{1}{4}R_{2 2} - \frac{1}{4}R_{3 1} - \frac{1}{4}R_{4}\\
	Y_5&=-\frac{1}{4}R_{1 1 1 1 1} + \frac{1}{4}R_{1 1 3} + \frac{1}{4}R_{1 2 2} + \frac{1}{4}R_{1 3 1} + \frac{1}{4}R_{2 1 2} + \frac{1}{4}R_{2 2 1} + \frac{1}{4}R_{3 1 1} - \frac{1}{4}R_{5}
\end{align}

Actually, the sums $F(n,k,r)$ are the coefficients of the ``amazing matrix'' of \cite{Ho}.
It is proved in this reference that
\begin{equation} \label{eq-fngen}
	F_b(n,k,r):= \sum_{i=0}^k(-1)^i{n+1\choose i}{n+bk-bi-r\choose n} = T_b(n+1, bk-r)
\end{equation}
where $T_b(n,p)$ is a $b$-nomial coefficient, {\it i.e.}, the coefficient of $x^p$
in $(1+x+\cdots+x^{b-1})^n$.

\section{Some transition matrices}

The construction of the basis $X_I$ is a special case of the following one.
Start with a multiplicative basis $\varphi^I$ of $\NCSF$, for example
$\varphi=\log\sigma_1$. Choose a formal series $f(z)=\sum_n c_nz^n$, and
define $Y:=\sum_nY_n=f(\varphi)$.
We have then
\begin{equation}\label{eq:Y2phi}
Y_n = \sum_{I\vDash n}c_{\ell(I)}\varphi^I
\end{equation}
so that

\begin{equation}
	Y^J = \sum_{I_1\vDash j_1,\ldots, I_r\vDash j_r}c_{\ell(I_1)}\cdots c_{\ell(I_r)}\varphi^{I_1\cdots I_r}.
\end{equation}
Let $\phi_I$ be the dual basis of $\varphi^I$ and $X_I$ be the dual basis
of $Y^I$. By duality, \eqref{eq:Y2phi} translates as

\begin{equation}\label{eq:phi2X}
	\phi_I = \sum_{I=I_1\cdots I_r}c_{\ell(I_1)}\cdots c_{\ell(I_r)}X_{|I_1|,\ldots,|I_r|}.
\end{equation}

By construction, the coproduct of $X_I$ is given by deconcatenation, since it is dual to a multiplicative basis,
so that the map $\Psi_f:\ \phi_I\mapsto X_I$ is a morphism of coalgebras. It is proved in \cite{FPT} that all
morphisms of coalgebras are actually of this type \cite[Th. 2.2]{FPT} and that more generally the group of coalgebra automorphisms
of a shuffle algebra is isomorphic to the group of invertible formal power series in one variable under composition
(in particular, $\Psi_f\circ\Psi_g=\Psi_{f\circ g}$).

This construction has been rediscovered in \cite{HI}, where the following convenient notation has been introduced.
If we encode compositions by words over an auxiliary alphabet $Z=\{z_i|i\ge 1\}$, endowed
with the operation $z_i\diamond z_j:=z_{i+j}$, and define a linear operator
$P$
by $P(z^I):=P_I$ for any basis $P_I$
of $Qsym$, \eqref{eq:phi2X} can be recoded as in \cite{HI}

\begin{equation}
	\phi\left(\frac{1}{1-\lambda z}\right)=
	X\left(\frac{1}{1-f_\diamond(\lambda z)}\right)
\end{equation}
where $z=\sum_i z_i$ and $f_\diamond(z)$ means that powers of $z$
are evaluated with the $\diamond$ operation\footnote{Actually, Hoffman and Ihara state this formula
for an arbitrary linear combination of the $z_i$, but this amounts only to a rescaling of the variables,
or for us, of the $\varphi_i$.}

If we choose $f$ such that $c_1=1$, we have $X_n=\phi_n$. We can then write

\begin{equation}
	\sum_I\lambda^{\ell(I)}\phi_I = \exp\left\{\lambda\sum_{k\ge 1}X_k\right\}.
\end{equation}
Indeed, the l.h.s. is the Cauchy kernel for the pair of bases $(\phi,\varphi)$ specialized
at the virtual alphabet $\Lambda$ defined by $\varphi_n(\Lambda)=\lambda$ for all $n$.
Then,

\begin{equation}
	\sigma_t(\Lambda)=\exp\left\{\lambda\sum_{k\ge 1}t^k\right\}=
	\exp\left\{\lambda\frac{t}{1-t}\right\},
\end{equation}
and
\begin{equation}
	\sigma_1(X\Lambda)=\exp\left\{\lambda\sum_{k\ge 1}\frac{x_i}{1-x_i}\right\}
	=\exp\left\{\lambda\sum_{k\ge 1}X_k\right\}.
\end{equation}
In the notation of \cite{HI,Kei}, this would be

\begin{equation}
	X(\exp_*(\lambda z))= \phi\left(\frac{1}{1-\lambda z}\right)
\end{equation}
where $*$ is the product rule for the $X$-basis, for example the quasi-shuffle
if $f(z)=e^z$ and the block shuffle of \cite{Kei} if $f(z)=\tanh z$.

\section{Generalisation}

In order to obtain interesting bases, the series $f(z)$ should have an addition
theorem of the type 
\begin{equation}
	f(x+y)=(f(x)+f(y))(1+\sum_{i+j\ge 1}c_{ij}f(x)^if(y)^j):= F(f(x),f(y)),
\end{equation}
({\it i.e.}, $F$ is a formal group law).

Tractable examples are rather scarce. To the list $e^z, \tanh z$ and $\tan z$,
one can add $f(z)=\frac 1q(e^{qz}-1)$, which satisfies
$f(x+y)=f(x)+f(y)+qf(x)f(y)$ and provides a deformation of the quasi-shuffle  \cite{FPT}
satisfying the recurrence
\begin{equation}
	au*bv = a(u*bv)+b(au*v)+q[a+b](u*v).
\end{equation}
These examples can be interpolated by the (known) formal group law
\begin{equation}\label{eq:FG}
	F_{\alpha,\beta}(x,y)=\frac{x+y+\alpha xy}{1-\beta xy}
\end{equation}
corresponding to the function
\begin{equation}
	f(x) = \frac{e^{sx}-e^{tx}}{se^{tx}-te^{sx}}
\end{equation}
where $\alpha=s+t$ and $\beta=st$.

This is the exponential generating function of the homogeneous Eulerian polynomials
\begin{equation}
f(x)=\sum_{n\ge 1}E_n(s,t)\frac{x^n}{n!}
\end{equation}
with
\begin{equation}
E_n(s,t)=\sum_{\sigma\in\SG_n}s^{r(\sigma)}t^{d(\sigma)}
\end{equation}
where $r(\sigma)$ is the number of rises of $\sigma$ and $d(\sigma)$ its number of descents.

\bigskip
Define $Y_n$ by $Y=f(\varphi)$. 
\begin{proposition}
The expansion of $Y_n$ on the $S$-basis is given by
\begin{equation}
	Y_n= \sum_{r=1}^n\frac{1}{r!}\left(\sum_{k=1}^r s(r,k)E_k(s,t)\right)S_n^{[r]}
\end{equation}
where $s(r,k)$ are the signed Stirling numbers of the first kind.
\end{proposition}

\Proof By definition,
\begin{equation}
	Y=\sum_{k\ge 1}E_k(s,t)\frac{\varphi^k}{k!}=\sum_{k\ge 1}E_k(s,t)E^{[k]}
\end{equation}
and since
\begin{equation}
	\sigma_1^x = \sum_kx^kE^{[k]} = \sum_{r\ge 0}{x\choose r}\sum_{\ell(I)=r}S^I
	=\sum_{r\ge 0}\frac1{r!}\sum_{k=0}^rs(r,k)x^k\sum_{\ell(I)=r}S^I
\end{equation}
we have
\begin{equation}
	E^{[k]}=\sum_r\frac1{r!}s(r,k)S^{[r]}
\end{equation}
whence the result.\qed

The expansion of the ribbon basis follows from the identity
\begin{equation}
\sum_{\ell(I)=r\atop |I|=n}S^I=\sum_{k=0}^{r-1}{n-r+k\choose k}{\A}(n,r-k).
\end{equation}

Another expression can be obtained as above in terms of the amazing matrix.
Imitating the calculation of Section \ref{sec:spec}, we have (treating $s-t$ as a scalar)

\begin{equation}
	Y_n=\frac1{s-t}\sum_{I\vDash n}\left(\frac{t}{s-t}\right)^{\ell(I)-1}S^I((s-t)A)
	= \sum_{J\vDash n}\frac{s^{n-\ell(J)}t^{\ell(J)-1}}{(s-t)^n}R_J((s-t)A)
\end{equation}
so that
\begin{equation}
	Y_n=\frac1{(s-t)^n}\sum_j s^{n-j}t^{j-1}\A(n,j)*S_n^{[s-t]}
\end{equation}
where
\begin{equation}
	\A(n,j)*S_n^{[s-t]}=\sum_iP_{ij}(s-t)\A(n,i).
\end{equation}

\begin{theorem}
Let  $X_I$ be the dual basis of $Y^I$.
	Then,
\begin{equation}
	X_I X_J = X(I\star J)
\end{equation}
where $\star$ is defined on words over the integers by the recursion
	\begin{equation}\label{eq:recab}
	au\star vb = a(u\star bv)+b(au\star v) +\alpha[a+b]\cdot u\star v + \beta\zeta_{a+b}(u\star v).
\end{equation}
In particular, the operation $\star$ is associative and commutative.
\end{theorem}

\Proof
By \eqref{eq:FG}, the coproduct of $Y$ is

\begin{equation}
\Delta Y = Y\otimes 1+1\otimes Y +\alpha\sum_{\ell(I)=\ell(J)}\beta^{\ell(I)-1}Y^I\otimes Y^J
	+ \sum_{|\ell(I)-\ell(J)|=1}\beta^{(\ell(I)+\ell(J)-1)/2}Y^I\otimes Y^J.
\end{equation}
Writing compositions as words to be in line with the notation of the theorem, we have
\begin{equation}
\<X_{au}X_{bv},Y^{cw}\> =\<X_{au}\otimes X_{bv},\Delta Y_c\Delta Y^{w}\>.
\end{equation}
If $c=a$, $\Delta Y_c$ contains $Y_a\otimes 1$ with coefficient 1, so that
\begin{align}
	\<X_{au}X_{bv},Y^{cw}\>&= \left\<X_{au}\otimes X_{bv},\left(Y_a\otimes 1+\cdots\right)\sum_{(w)}Y^{w_{(1)}}\otimes Y^{w_{(2)}}\right\>\\
	&=\sum_{(w)}\<X_{au},Y_aY^{w_{(1)}}\>\<X_{bv},Y^{w_{(2)}}\>\\
	&=\sum_{(w)}\<X_{u},Y^{w_{(1)}}\>\<X_{bv},Y^{w_{(2)}}\>\\
	&=\sum_{(w)}\<X_{u}\otimes X_{bv},Y^{w_{(1)}}\otimes Y^{w_{(2)}}\>\\
	&= \<X_uX_{bv},Y^w\>.
\end{align}
Symmetrically, if $c=b$, $\<X_{au}X_{bv},Y^{cw}\>=\<X_{au}X_{v},Y^w\>$.

If $c=a+b$, $Y_a\otimes Y_b$ occurs in $\Delta Y_c$ with coefficient $\alpha$, and a similar calculation
shows that  $\<X_{au}X_{bv},Y^{cw}\>=\alpha \<X_uX_v,Y^w\>$.

If $c=a+b+d>a+b$, then the sum of the terms starting with $Y_a\otimes Y_b$ in $\Delta Y^c$ is
\begin{equation}
	\sum_{|\ell(s)-\ell(t)|=1\atop |s|+|t|=d}Y^{as}\otimes Y^{bt}= (Y_a\otimes Y_b)\Delta Y_d
\end{equation}
which occurs with  coefficient $\beta$. Thus,
\begin{align}
	\<X_{au}X_{bv},Y^{cw}\>&=\<Y_{au}\otimes Y_{bv}, (Y_a\otimes Y_b)\Delta Y_d \Delta Y^w\>\\
	&=\<X_u\otimes X_v, \Delta(Y_dY^w)\>\\
	&= \<X_uX_v, Y_dY^w\>\\
	&= \<X_uX_v,\partial_{a+b}Y^{cw}\>
\end{align}
where the map  $\partial_xY^{a_1a_2\cdots a_r}:=Y^{(a_1-x)\cdot a_2\cdots a_r}$ is the adjoint of $\zeta_x$,
so that finally
\begin{equation}
	\<X_{au}X_{bv},Y^{cw}\>=\<\zeta_{a+b}(X_uX_v), Y^{cw}\>.
\end{equation}
In all other cases, it is clear that  $\<X_{au}X_{bv},Y^{cw}\>=0$, whence the recurrence \eqref{eq:recab}.
\qed

\bigskip
Clearly, $X_1=\phi_1=F_1=M_1$ is the sum of the variables. It is always interesting to expand its powers
on a new basis. Since $\phi_{1^n}=\phi_1^n$, it follows from \eqref{eq:phi2X} that

\begin{equation}
	X_1^n = \sum_{I \vDash n}{n\choose I} E^I(s,t) X_I
\end{equation}
where for $I=(i_1,\ldots,i_r)$, ${n\choose I}={n\choose i_1\ldots i_r}$ is the multinomial coefficient
and $E^I=E_{i_1}\cdots E_{i_r}$ is a product of Eulerian polynomials.

\bigskip
For $I=1^n$, $X_{1^n}$ is a symmetric function. The coefficients of its monomial expansion are

\begin{equation}
	\<X_{1^n},S^I\>=\<\Delta X_{1^n},S_{i_1}\otimes\cdots\otimes S_{i_r}\>=\prod_{k=1}^r\<X_{1^{i_k}},S_{i_k}\>=\prod_{k=1}^rd_{i_k}
\end{equation}
where

\begin{equation}
	d_n = \<X_{1^n},S_n\> =[y^n]\left(\frac{1+sy}{1+ty}\right)^\frac{1}{s-t}.
\end{equation}

\bigskip
One can give a closed formula for the product $X_{1^p}X_{1^q}$. 
\begin{proposition}
	For a composition $I$, let $\ell_0(I)$ denote its number of even parts,
	and $\ell_1(I)$ its number of odd parts. Then
\begin{equation}
	X_{1^p}X_{1^q} = \sum_{I\vDash p+q}{\ell_1(I)\choose p-\frac{p+q-\ell_1(I)}{2}}\alpha^{\ell_0(I)}\beta^{\frac{p+q-\ell_1(I)}{2}}X_I.
\end{equation}
\end{proposition}
\Proof
By duality,
\begin{align}
	\<X_{1^p}X_{1^q},Y^I\>&= \<\Delta^r(X_{1^p}X_{1^q}), Y_{i_1}\otimes\cdots\otimes Y_{i_r}\>\\
	&= \sum_{p_1+\cdots+p_r=p,q_1+\cdots+q_r=q\atop p_1+q_1=i_1,\ldots, p_r+q_r=i_r}\prod_{k=1}^r\<X_{1^{p_k}}X_{1^{q_k}},Y_{i_k}\>.
\end{align}
Now,
\begin{align}
	\<X_{1^p}X_{1^q},Y_n\>&=\<X_{1^p}\otimes X_{1^q},\Delta Y_n\>\\
	&= \left\<X_{1^p}\otimes X_{1^q}, \frac{Y\otimes 1+1\otimes Y+\alpha Y\otimes Y}{1\otimes 1-\beta Y\otimes Y}\right>\\
	&= \left\<X_{1^p}\otimes X_{1^q}, \frac{Y_1\otimes 1+1\otimes Y_1+\alpha Y_1\otimes Y_1}{1\otimes 1-\beta Y_1\otimes Y_1}\right>\\
	&=\begin{cases}
		\alpha \beta^{\frac{n}2-1}&\text{if $p=q$, so that $n=p+q$ is even,}\\
		\beta^{\frac{n-1}{2}}&\text{if $|p-q|=1$, so that $n=p+q$ is odd.}
	\end{cases}
\end{align}
Thus, the coefficient of $X_I$ is $X_{1^p}X_{1^q}$ is equal to the number of $2\times r$ integer matrices
\begin{equation}
	\left[\begin{matrix}p_1&p_2&\cdots&p_r\\
		            q_1&q_2&\cdots&q_r
	\end{matrix}\right]
\end{equation}
with row sums $p,q$, column sums $i_1,\ldots,i_r$, and such that $p_k=q_k$ if $i_k$ is even,
and $|p_k-q_k|=1$ of $i_k$ is odd. On can form such a matrix by adding to the matrix $p'_k=q'_k=\lfloor\frac{i_k}{2}\rfloor$
a matrix of 0 and 1, with 0 in the columns of the even $i_k$, and exactly one 1 in the columns of the odd $i_k$, such that
the sum of the first row is $p-\frac{p+q-\ell_1(I)}{2}$ and, equivalently, such that of the second row is
 $q-\frac{p+q-\ell_1(I)}{2}$, whence the binomial coefficient. \qed

\section{Extension to other algebras}

The algebra of quasi-symmetric functions has a noncommutative version $\WQSym$  (Word Quasi-Symmetric functions)
consisting of the invariants of Hivert's quasi-symmetrizing action of the symmetric group on the free associative algebra \cite{NCSF7}.
The \emph{packed word} $u=\pack(w)$ associated with a word $w\in A^*$ is
obtained by the following process. If $b_1<b_2<\ldots <b_r$ are the letters
occuring in $w$, $u$ is the image of $w$ by the homomorphism
$b_i\mapsto a_i$.
A word $u$ is said to be \emph{packed} if $\pack(u)=u$. We denote by $\PW$ the
set of packed words.
With such a word, we associate the polynomial
\begin{equation}
\M_u :=\sum_{\pack(w)=u}w\,.
\end{equation}
The product of the monomial basis of $\WQSym$ is described by convolution of packed words
\begin{equation}
	\M_u\M_v=\sum_{w=u'v'\atop \pack(u')=u,\ \pack(v')=v}\M_w
\end{equation}
which can also be described in terms of the bigger algebra $\MQSym$ (Matrix Quasi-Symmetric functions) \cite{NCSF6}, based on
packed integer matrices. The product of two basis elements $\MS_P\MS_Q$, where $P$ is a $p\times r$ matrix and
$Q$ a $q\times s$ matrix is obtained by forming the block matrix
\begin{equation}
	P\bullet Q =\left[\begin{matrix}P&0\\0&Q\end{matrix}\right] 
\end{equation}
regarded as a word over the alphabet of rows, and taking the quasi-shuffle of the word $P_1\cdots P_r$ formed by its
first $p$ rows with the word formed by its last $q$ rows, the contractions being given by vector addition of the rows.

Packed words can be encoded by packed $0-1$-matrices with exactly one 1 in each column. For example, $u=21321$ is encoded by
\begin{equation}
	\left[\begin{matrix}0&1&0&0&1\\
	                    1&0&0&1&0\\
	                    0&0&1&0&0 
	\end{matrix}\right]
\end{equation}
where $m_{ij}=1$ iff the $j$th letter is an $i$. Then, packed convolution corresponds to the quasi-shuffle of these matrices,
and so can be described by a recurrence 
\begin{equation}
	AP\stuffle BQ = A'(P\stuffle BQ)+B''(AU\stuffle Q) + (A'+B'')(P\stuffle Q)
\end{equation}
where $A'$ and $B''$ denote $A$ and $B$ completed by the appropriate number of zeros on the right or on the left.

The basis $X_I$ can be lifted to $\WQSym$ by means of the (reverse) refinement order on packed words: if 
\begin{equation}
	X_I=\sum_{J\le I}c_{IJ}M_J\quad\text{we define}\quad \XX_u=\sum_{v\le u}c_{\ev(u)\ev(v)}\M_v.
\end{equation}

To describe the product $\XX_u\XX_v$, we shall embed $\WQSym$ and $\MQSym$ in the limit $\ell\rightarrow\infty$
of the level $\ell$ quasi-symmetric functions defined in \cite{NTcol}, where the noncommutative product of these
algebras can be realized as a shifted version of the commutative product of a bigger algebra.

\def\bi{{\bf i}}
\def\bj{{\bf j}}
\def\bk{{\bf k}}
\def\bn{{\bf n}}
\def\bI{{\bf I}}
\def\bJ{{\bf J}}
\def\bK{{\bf K}}

Let $\Omega$ be the semigroup of nonnegative integer sequences $\bi =(i_0,i_1,\ldots)$ with finite sum
$|\bi|$. We denote by $\Sym^{(\N)}$ the Hopf algebra of $\N$-colored noncommutative symmetric functions,
which is defined as the free associative algebra over indeterminates $S_{\bi}$, with $S_{\bf 0}=1$,
endowed with the coproduct
\begin{equation}
	\Delta S_\bn = \sum_{\bi+\bj=\bn}S_\bi\otimes S_\bj.
\end{equation}
Departing from the notation of \cite{NTcol}, we represent color sequences $\bi$ as row vectors, and regard the label $\bI$ of
the basis element $S^\bI=S_{\bi_1}\cdots S_{\bi_r}$ as an $r\times\infty$ matrix. The multidegree $||\bI||$ of $\bI$
is defined as its column sum sequence. The Hopf algebra $QSym^{(\N)}$ of $\N$-colored quasi-symmetric functions
is defined as its graded dual with respect to this multidegree.

Both algebras admit polynomial realizations, in terms of two colored alphabets

\begin{equation}
	\A = \bigsqcup_{c\ge 0}A^{(c)}, \quad \X = \bigsqcup_{c\ge 0}X^{(c)},
\end{equation}

In terms of $\A$, the generating function of the complete functions can be
written as
\begin{equation}
\sigma_{\bf x}(\A) = \prod_{i\geq1}^{\rightarrow}
\left(1-\sum_{j\ge 0} x^{(j)} a_{i}^{(j)} \right)^{-1}
= \sum_\npn {S_{\bf n}(\A) {\bf x}^{\bf n}},
\end{equation}
where ${\bf x}^{\bf n} = (x^{(0)})^{n_0} \cdots (x^{(k)})^{n_{k}}\cdots $.

This realization gives rise to a Cauchy formula 
which in turn allows one to identify the dual of $\NCSF^{(\N)}$
with the limit of an algebra introduced by S. Poirier in~\cite{Poi}.

Let $\X = X^{(0)} \sqcup \cdots \sqcup X^{(k)}\sqcup\cdots$,
where $X^{(i)}=\{ x_j^{(i)},j\geq1\}$, be an $\N$-colored alphabet of
commutative variables, also commuting with $\A$.
Imitating the level $1$ case (see~\cite{NCSF6}), we define the Cauchy kernel

\begin{equation}
\label{CauchyXA}
K(\X,\A) = \prod_{j\geq1}^{\rightarrow}
\sigma_{\left(x_j^{(0)}, \ldots, x_j^{(k)},\ldots \right)} (\A).
\end{equation}

Expanding on the basis $S^{\bf I}$ of $\NCSF^{(\N)}$, we get as coefficients
what can be called the \emph{$\N$-monomial quasi-symmetric functions}
$M_{\bf I}(\X)$

\begin{equation}
\label{K-expand}
K(\X,\A) = \sum_{\bf I} M_{\bf I}(\X) S^{\bf I}(\A),
\end{equation}
defined by
\begin{equation}
M_{\bf I}(\X) = \sum_{j_1<\cdots<j_m}
{\bf x}^{{\bf i}_1}_{j_1} \cdots
{\bf x}^{{\bf i}_m}_{j_m},
\end{equation}
with ${\bf I}=({\bf i}_1,\ldots,{\bf i}_m)$.

These functions form a basis of a subalgebra $\QSym^{(\N)}$ of
$\K[\X]$, which we shall call the \emph{algebra of $\N$-colored quasi-symmetric functions }.

We can now define $\varphi_{\bf x}=\log\sigma_{\bf x}$, and get a family of primitive
generators of $\Sym^{(\N)}$ by setting $\varphi_\bn =$ coefficient of ${\bf x}^\bn$ in
$\varphi_{\bf x}$. 

With a formal series $f(x)=x+O(x^2)$ as above, we can now define a basis $Y^\bI$
by setting $Y_\bn=$ coefficient of ${\bf x}^\bn$ in $f(\varphi_{\bf x})$. The coproduct
of $Y$ will the be given by the formal group law associated with $f$.

Thus, if $f$ is as above chosen as  the exponential generating functions of the $E_n(s,t)$,
the product rule of basis $X_\bI$ will be the $(\alpha,\beta)$-quasi-shuffle over the alphabet
$\Omega$.

The subspace of $QSym^{(\N)}$ spanned by the $M_\bI$ where $\bI$ is a packed matrix (followed by zero columns) is stable
for the product. This is not the case of the span of the matrices encoding packed words,
but both are stable for the shifted product: if $\bI$ is an $r\times p$ packed matrix, and $\bJ$ an $s\times q$
packed matrix, define
\begin{equation}
	M_\bI\star M_\bJ = M_\bI\cdot M_{{\bf 0}_s^p\bJ}
\end{equation}
that is, shift $\bJ$ by the appropriate number of zero columns so as to reproduce the product
of $\MQSym$. Then, since $X_\bI$ is a linear combination of $M_\bJ$ involving only
coarser $\bJ$ (obtained from $\bI$ by adding consecutive rows), we have as well
\begin{equation}
	X_\bI\star X_\bJ = X_\bI\cdot X_{{\bf 0}_s^p\bJ}
\end{equation}
This is therefore an encoding of the product $\XX_u\XX_v$ in $\WQSym$, which is thus
given by the shifted $(\alpha,\beta)$-quasi-shuffle of the corresponding matrices.
Since in this case no two contractions can be equal, the product is multiplicity-free,
and the cofficient of each $\XX_w$ is just a single monomial $\alpha^i\beta^j$, which
can be explicitly computed.

Let $f_{uv}^w$ be this coefficient. If there exists $j\not=k$ such that
$u_j=u_k$ and $w_j\not=w_k$, or $v_j=v_k$ and $w_{j+|u|}\not=w_{k+|u|}$,
then $f_{uv}^w=0$. Otherwise, let $a_i$ be the number of different values
$u_j$ such that $w_j=i$, and $b_i$ be the number of different values $v_j$
such that $w={j+|u|}=i$. Then,

\begin{equation}
	f_{uv}^w = \prod_{i=1}^{\max(w)}f_{uv}^w(i),\quad\text{where}\ 
	f_{uv}^w(i)=
	\begin{cases}
		0&\text{if $|a_i-b_i|>1$,}\\
		\alpha\beta^{a_i-1}&\text{if $a_i=b_i$,}\\
		\beta^{\min(a_i,b_i)}&\text{otherwise}.
	\end{cases}
\end{equation}
For example, the coefficient of $\XX_{11211}$ in $\XX_{123}\XX_{11}$ is $\alpha\beta$,
since $a_1=2,b_1=2$ and $a_2=1,b_2=0$.

\section{Appendix: notations and background}\label{sec:app}
The graded dual of the Hopf algebra $QSym$ is the Hopf algebra $\NCSF$ of noncommutative symmetric functions \cite{NCSF1}. As a graded algebra, its is freely generated by a sequence $(S_n)_{n\ge 1}$ of noncommuting variables ($S_n$ being of degree $n$), called complete homogenous functions. They can be realized by means of an auxiliary alphabet $A=\{a_i|i\ge 1\}$ of noncommuting variables, and setting
\begin{equation}
\sigma_t(A) := \sum_{n\ge e 0}t^n S_n(A) = \prod^\rightarrow_{i\ge 1}(1-ta_i)^{-1}.
\end{equation}
The dual basis of the quasi-monomial basis $M_I$ of $QSym$ is the basis $S^I=S_{i_1}\cdots S_{i_r}$
of $\NCSF$. In general, we note with a superscript such {\it multplicative bases}, formed with product of algebraic generators. Other examples are $\Lambda^I$ and $\Phi^I$, where the elementary functions $\Lambda_n$ and the $\Phi_n$ are defined by the series
\begin{equation}
\lambda_t(A)=\sum_{n\ge 0}t^n\Lambda_n(A)=\prod^\leftarrow_{i\ge 1}(1+ta_i)
\end{equation}
and
\begin{equation}
\varphi(t) = \log \sigma_t = \sum_{n\ge _1}t^n\frac{\Phi_n}{n}.
\end{equation}
The dual basis of the fundamental basis $F_I$ of $QSym$ is the ribbon basis $R_I$ of $\NCSF$.
It is realized as the sum of words with descent composition $I$.

The internal product $*$ of $\NCSF$ is dual to the internal coproduct $\delta F(X)=F(XY)$ of $QSym$,
and $(\NCSF_n,*)$ is anti-isomorphic to Solomon's descent algebra $\Sigma_n$, $R_I$ being interpreted as the sum of permutations with descent composition $I$.

For a scalar $\alpha$, the symmetric functions of the virtual alphabet $\alpha A$ are defined by
$\sigma_t(\alpha A)=\sigma_t(A)^\alpha$. 

The {\em Eulerian idempotents} $E_n^{[k]}$ are the homogenous components
of degree $n$ in the series $E^{[k]}$ defined by
\begin{equation}
\sigma_t(A)^x=\sum_{k\ge 0}x^k E^{[k]}(A),
\end{equation}
(see \cite[Section 5.3]{NCSF1}).
We have
\begin{equation}
E_n^{[k]}*E_n^{[l]} = \delta_{kl}E_n^{[k]}\,,\quad\text{and}\quad\sum_{k=1}^nE_n^{[k]}=S_n,
\end{equation}
so that the $E_n^{[k]}$ span a commutative $n$-dimensional $*$-subalgebra of
$\Sym_n$, called the Eulerian subalgebra.

The {\em noncommutative Eulerian polynomials} are defined by \cite[Section 5.4]{NCSF1}
\begin{equation}
{\mathcal A}_n(t) =
\ \sum_{k=1}^n \ t^k\, \Big(
\sum_{{\scriptstyle |I|=n}\atop{\scriptstyle \ell(I)=k}} R_I \,
\Big)
=
\ \sum_{k=1}^n \ {\A}(n,k)\, t^k \, .
\end{equation}
The generating series of the ${\mathcal A}_n(t)$ is
\begin{equation}
{\mathcal A}(t) := \ \sum_{n\ge 0} \, {\mathcal A}_n(t)
=
(1-t) \, \left( 1 - t\, \sigma_{1-t} \right)^{-1} \,,
\end{equation}
where $\sigma_{1-t}=\sum (1-t)^nS_n$.

\medskip
Let ${\mathcal A}_n^*(t) = (1-t)^{-n}\, {\mathcal A}_n(t)$.
Then,
\begin{equation}
{\mathcal A}^*(t)
:=
\ \sum_{n\ge 0} \, {\mathcal A}_n^*(t)
=
\sum_{I} \
\left( \displaystyle {t \over 1-t} \right)^{\ell(I)} \, S^I \ .
\end{equation}
This last formula can also be written in the form
\begin{equation} \label{GEN*}
{\mathcal A}^*(t)
=
\ \sum_{k\ge 0} \ \left(
{t\over 1-t}\right)^k \left( S_1+S_2+S_3+\cdots\, \right)^k
\end{equation}
or
\begin{equation}\label{GEN_A}
{1\over 1-t\, \sigma_1(A)}
=
\ \sum_{n\ge 0}\ {{\mathcal A}_n(t)\over (1-t)^{n+1}} \ .
\end{equation}
Let $S^{[k]}=\sigma_1(A)^k$ be the coefficient of $t^k$ in this series. In degree $n$,
\begin{equation}\label{S2E}
S_n^{[k]}=\sum_{I\vDash n, \ell(I)\le k}{k\choose \ell(I)}S^I=\sum_{i=1}^nk^iE_n^{[i]}\,.
\end{equation}
Conversely,
\begin{equation}
{ {\mathcal A}_n(t) \over (1-t)^{n+1} }
=
\ \sum_{k\ge 0}\ t^k\, S_n^{[k]} \ ,
\end{equation}
so that
\begin{equation}\label{A2S}
{\A}(n,p)
=
\ \sum_{i=0}^p\ (-1)^i\, {n+1\choose i}\, S_n^{[p-i]} \ .
\end{equation}
%



\footnotesize
\begin{thebibliography}{aa}
%

\bibitem{NCSF6} {\sc G. Duchamp, F. Hivert}, and {\sc J.-Y. Thibon}, 
{\it Noncommutative symmetric functions VI: free quasi-symmetric functions and related algebras},
Internat. J. Alg. Comput. {\bf 12} (2002), 671--717.

\bibitem{NCSF7}{\sc G. Duchamp, F. Hivert, J.-C. Novelli} and {\sc J.-Y. Thibon}, {\it
	Noncommutative symmetric functions VII: free quasi-symmetric functions revisited},
		Ann. Comb. {\sc 15} (2011), no. 4, 655--673.

\bibitem{FPT}{\sc L. Foissy, F. Patras} and {\sc J.-Y. Thibon}, {\it
Deformations of shuffles and quasi-shuffles}, Ann. Institut Fourier {\bf 66} (2016), 209--327.

\bibitem{NCSF1} {\sc I.~M. Gelfand, D. Krob, A. Lascoux, B. Leclerc,
V.~S. Retakh}, and {\sc J.-Y. Thibon}.
{\it Noncommutative symmetric functions},
Adv. Math. {\bf 112}, 1995, 218--348.
%
%
\bibitem{Ge} {\sc I.~Gessel}, {\it Multipartite P-partitions and inner
products of skew Schur functions}, [in ``Combinatorics and algebra",
C. Greene, Ed.], Contemporary Mathematics, {\bf 34} (1984), 289--301.
%

\bibitem{HHL}{\sc 
	J. Haglund, M. Haiman} and {\sc  N. Loehr},
		{\it A combinatorial formula for Macdonald polynomials}, J. Am. Math. Soc. {\bf 18} (2005) 735--761.

\bibitem{HNTT}{\sc F. Hivert, J.-C. Novelli, L. Tevlin} and {\sc J.-Y. Thibon}, {\it
Permutation statistics related to a class of noncommutative symmetric functions and generalizations of the Genocchi numbers},
		Selecta Math. (N.S.) {\bf 15} (2009), no. 1, 105--119.

\bibitem{Hof}{\sc M. Hoffman}, {\it
	Quasi-shuffle products}, J. Algebraic Combin. {\bf 11} (2000), no. 1, 49--68. 

\bibitem{HI}{\sc M. Hoffman} and {\sc K. Ihara}, {\it
	Quasi-shuffle products revisited}, J. Algebra {\bf 481} (2017), 293--326.

\bibitem{Ho}{\sc J. Holte}, {\it Carries, combinatorics, and an amazing matrix}, 
Amer. Math. Mon.  {\bf 104} (1997) 138--149.


\bibitem{Kei} {\sc A. Keilthy}, 
	{\it A generalisation of quasi-shuffle algebras and an application
		to multiple zeta values}, INTEGERS 22 A114 (2022);  arXiv:2202.04739.
%

\bibitem{NCSF2} {\sc Krob D., Leclerc B.}, and {\sc Thibon J.-Y.},
{\it Noncommutative symmetric functions II: Trans\-for\-ma\-tions of
alphabets}, Intern. J. Alg. Comput.,
{\bf 7}, (2), (1997), 181--264.
%
%
\bibitem{NCSF4} {\sc D.~Krob} and {\sc J.-Y.Thibon},
{\it Noncommutative symmetric functions IV{\,}: 
  Quantum linear groups and Hecke algebras at $q=0$}, 
J. Alg. Comb. {\bf 6} (1997), 339--376.
%
\bibitem{LNT}{\sc A. Lascoux, J.-C. Novelli} and {\sc J.-Y. Thibon}, {\it
	Noncommutative symmetric functions with matrix parameters}, J. Algebraic Combin. {\bf 37} (2013), no. 4, 621--642.
%
\bibitem{Mac} {\sc Macdonald I.G.}, {\it Symmetric functions and Hall
polynomials}, Clarendon Press, 1995.
%

\bibitem{Mason}{\sc S. K. Mason}, {\it 
Recent trends in quasisymmetric functions},
Recent trends in algebraic combinatorics, 239--279, Assoc. Women Math. Ser., {\bf 16}, Springer, Cham, 2019.

\bibitem{NPT}{\sc J.-C. Novelli, F. Patras} and {\sc J.-Y. Thibon}, {\it
	Natural endomorphisms of quasi-shuffle Hopf algebras}a, Bull. Soc. Math. France {\bf 141} (2013), no. 1, 107--130. 

\bibitem{NTbinbsh}{\sc J.-C. Novelli} and {\sc J.-Y. Thibon},
	{\it Binary shuffle bases for quasi-symmetric functions},
Ramanujan journal {\bf 40} (2016), 207--225. 
%

\bibitem{NTcol}{\sc J.-C. Novelli} and {\sc J.-Y. Thibon},{\it
	Free quasi-symmetric functions and descent algebras for wreath products, and noncommutative multi-symmetric functions},
	Discrete Math. {\bf 310} (2010), no. 24, 3584--3606. 



\bibitem{NTLag2} {\sc J.-C. Novelli}  and {\sc J.-Y. Thibon},{\it 
Noncommutative symmetric functions and Lagrange inversion II: noncrossing partitions and the Farahat-Higman algebra},
		Adv. in Appl. Math. {\bf 140} (2022), Paper No. 102396, 39 pp.

\bibitem{NTamaz}{\sc J.-C. Novelli} and {\sc J.-Y. Thibon},{\it
	Noncommutative symmetric functions and an amazing matrix}, Adv. in Applied Math. {\bf 48} (2012) 528--534. 

\bibitem{NTT}{\sc J.-C. Novelli, L. Tevlin} and {\sc J.-Y. Thibon}, {\it
On some noncommutative symmetric functions analogous to Hall-Littlewood and Macdonald polynomials},
Internat. J. Algebra Comput. {\bf 23} (2013), no. 4, 779--801.

\bibitem{NTToum} {\sc J.-C. Novelli, J.-Y. Thibon} and {\sc F. Toumazet}, {\it
A noncommutative cycle index and new bases of quasi-symmetric functions and noncommutative symmetric functions},
Ann. Comb. {\bf 24} (2020), no. 3, 557-576.

\bibitem{NTW}{\sc J.-C. Novelli, J.-Y. Thibon} and {\sc L. K. Williams}, {\it
	Combinatorial Hopf algebras, noncommutative Hall-Littlewood functions, and permutation tableaux},
Adv. Math. {\bf 224} (2010), no. 4, 1311--1348. 

\bibitem{Poi}{\sc S. Poirier}, {\it
	Cycle type and descent set in wreath products},
	Proceedings of the 7th Conference on Formal Power Series and Algebraic Combinatorics (Noisy-le-Grand, 1995). Discrete Math. 180 (1998), no. 1-3, 315--343. 

\bibitem{SchaThi}{\sc T. Scharf} and {\sc J.-Y. Thibon}, {\it 
	On Witt vectors and symmetric functions}, Algebra Colloq. {\bf 3} (1996), no. 3, 231--238.


%


\bibitem{TU}{\sc J.-Y.~Thibon} and {\sc B.-C.-V.~Ung}, 
{\it Quantum quasi-symmetric functions and Hecke algebras}, 
J. Phys. A: Math. Gen., {\bf 29} (1996),
7337--7348.

\bibitem{Zag}{\sc D. Zagier},{\it
Values of zeta functions and their applications}, 
in `First European Congress of Mathematics'' Vol. II, pp. 49--512, Birkhauser Boston, Cambridge, MA, 1994.
\end{thebibliography}
\end{document}